\documentclass{wscpaperproc}
\usepackage{latexsym}
\usepackage{graphicx}
\usepackage{mathptmx}
\usepackage{float}
\usepackage{tikz}
\usepackage{pgfplots}
\usepackage{pgf}
\usepackage{tikz-inet}

\usepackage{amsmath}
\usepackage{amsfonts}
\usepackage{amssymb}
\usepackage{amsbsy}
\usepackage{amsthm}
\usepackage{tasks}

\usepackage[pdftex,colorlinks=true,urlcolor=blue,citecolor=black,anchorcolor=black,linkcolor=black]{hyperref}

\usepackage{algorithm}
\usepackage[noend]{algpseudocode}
\usepackage{array,multirow,makecell}
\usepackage{graphicx} 

\newtheoremstyle{wsc}
{3pt}
{3pt}
{}
{}
{\bf}
{}
{.5em}
{}

\theoremstyle{wsc}

\def\cO {{\mathcal{O}}}
\def\cP {{\mathcal{P}}}

\def\cX {{\mathcal{X}}}

\def\bU {{\mathbf{U}}}
\def\bV {{\mathbf{V}}}

\def\EE {{\mathbb E}}

\def\II {{\mathbb I}}

\def\RR {{\mathbb R}}

\def\rma{{\rm a}}

\def\d{{\rm d}}
\def\rme{{\rm e}}
\def\rmy{{\rm y}}
\def\arqmc{{\text{arqmc}}}
\def\mua{\hat{\mu}_{n}^{\text{arqmc}}}
\def\Var{{\rm Var }}
\def\arqmc{{\rm arqmc}}

\def\q {$\kern15pt$}    
\def\?{\discretionary{}{}{}}  

\def\g{\leftarrow}   
\def\dGamma{{\rm Gamma}}
\def\dNormal{{\rm Normal}}
\def\dUnif{{\rm Uniform}}

\def\hfuture#1{}

\newif\ifnotes\notestrue

\def\hpierre#1{}
\newcommand{\hmpierre}[1]{{}}

\def\hamal#1{}

\def\hflorian#1{}

\pgfplotscreateplotcyclelist{defaultcolorlist}{%
	{blue!95!black,line width=0.9pt,mark=dot*,solid},
	{red!95!black,line width=0.9pt,mark=square*,solid},
	{green!98!black,line width=0.9pt,mark=triangle*,solid},
	{black,mark=star*,solid},
	{brown!85!black,mark=square,dashed},
	{purple!85!black,mark=triangle,dotted}}

\begin{document}

	%
	%

\pagestyle{fancyplain}

\thispagestyle{plain}
\firstPageHead{}

\chead{\fancyplain{}{\itshape Ben Abdellah, L'Ecuyer, and Puchhammer}}

\rhead{}
\cfoot{}
\renewcommand{\headrulewidth}{0pt} 

\makeatletter
\let\@internalcite\cite
\def\cite{\def\@citeseppen{-1000}%
    \def\@cite##1##2{(##1\if@tempswa , ##2\fi)}%
    \def\citeauthoryear##1##2##3{##1 ##3}\@internalcite}
\def\citeNP{\def\@citeseppen{-1000}%
    \def\@cite##1##2{##1\if@tempswa , ##2\fi}%
    \def\citeauthoryear##1##2##3{##1 ##3}\@internalcite}
\def\citeN{\def\@citeseppen{-1000}%
    \def\@cite##1##2{##1\if@tempswa, ##2)\else{}\fi}%
    \def\citeauthoryear##1##2##3{##1 (##3)}\@citedata}
\def\citeA{\def\@citeseppen{-1000}%
    \def\@cite##1##2{(##1\if@tempswa , ##2\fi)}%
    \def\citeauthoryear##1##2##3{##1}\@internalcite}
\def\citeANP{\def\@citeseppen{-1000}%
    \def\@cite##1##2{##1\if@tempswa , ##2\fi}%
    \def\citeauthoryear##1##2##3{##1}\@internalcite}
\def\shortcite{\def\@citeseppen{-1000}%
    \def\@cite##1##2{(##1\if@tempswa , ##2\fi)}%
    \def\citeauthoryear##1##2##3{##2 ##3}\@internalcite}
\def\shortciteNP{\def\@citeseppen{-1000}%
    \def\@cite##1##2{##1\if@tempswa , ##2\fi}%
    \def\citeauthoryear##1##2##3{##2 ##3}\@internalcite}
\def\shortciteN{\def\@citeseppen{-1000}%
    \def\@cite##1##2{##1\if@tempswa, ##2\else{}\fi}%
    \def\citeauthoryear##1##2##3{##2 (##3)}\@citedata}
\def\shortciteA{\def\@citeseppen{-1000}%
    \def\@cite##1##2{(##1\if@tempswa , ##2\fi)}%
    \def\citeauthoryear##1##2##3{##2}\@internalcite}
\def\shortciteANP{\def\@citeseppen{-1000}%
    \def\@cite##1##2{##1\if@tempswa , ##2\fi}%
    \def\citeauthoryear##1##2##3{##2}\@internalcite}
\def\citeyear{\def\@citeseppen{-1000}%
    \def\@cite##1##2{(##1\if@tempswa , ##2\fi)}%
    \def\citeauthoryear##1##2##3{##3}\@citedata}
\def\citeyearNP{\def\@citeseppen{-1000}%
    \def\@cite##1##2{##1\if@tempswa , ##2\fi}%
    \def\citeauthoryear##1##2##3{##3}\@citedata}
%
%
%
\def\@citedata{%
    \@ifnextchar [{\@tempswatrue\@citedatax}%
                  {\@tempswafalse\@citedatax[]}%
}

\def\@citedatax[#1]#2{%
\if@filesw\immediate\write\@auxout{\string\citation{#2}}\fi%
  \def\@citea{}\@cite{\@for\@citeb:=#2\do%
    {\@citea\def\@citea{, }\@ifundefined
       {b@\@citeb}{{\bf ?}%
       \@warning{Citation `\@citeb' on page \thepage \space undefined}}%
{\csname b@\@citeb\endcsname}}}{#1}}%

%
\def\@citex[#1]#2{%
\if@filesw\immediate\write\@auxout{\string\citation{#2}}\fi%
  \def\@citea{}\@cite{\@for\@citeb:=#2\do%
    {\@citea\def\@citea{, }\@ifundefined
       {b@\@citeb}{{\bf ?}%
       \@warning{Citation `\@citeb' on page \thepage \space undefined}}%
{\csname b@\@citeb\endcsname}}}{#1}}%

%
\def\@biblabel#1{}
\makeatother



\newdimen\bibindent
\bibindent=0.0em
\def\thebibliography#1{\section*{\refname}\list
   {}{\settowidth\labelwidth{[#1]}
   \leftmargin\parindent
   \itemindent -\parindent
   \listparindent \itemindent
   \itemsep 0pt
   \parsep 0pt}
   \def\newblock{}
   \sloppy
   \sfcode`\.=1000\relax}


\setlength{\baselineskip}{12.7pt}

	\title{ ARRAY-RQMC FOR OPTION PRICING UNDER STOCHASTIC VOLATILY MODELS }
	
	\author{Amal Ben Abdellah\\ 
		Pierre L'Ecuyer \\ 
		Florian Puchhammer\\ \\ 
		D\'epartement d'Informatique et de Recherche Op\' erationnelle \\
		Pavillon Aisenstadt, Universit\'e de Montr\'eal, C.P.6128, Succ. Centre-Ville \\
		Montr\'eal (Qu\'ebec), H3C 3J7, CANADA\\
	}
	
	\maketitle
	
\section*{ABSTRACT}

Array-RQMC has been proposed as a way to effectively apply randomized quasi-Monte Carlo (RQMC) when simulating a Markov chain over a large number of steps to estimate an expected cost or reward. The method can be very effective when the state of the chain has low dimension. For pricing an Asian option under an ordinary geometric Brownian motion model, for example, Array-RQMC reduces the variance by huge factors. In this paper, we show how to apply this method and we study its effectiveness in case the underlying process has stochastic volatility. We show that Array-RQMC can also work very well for these models, even if it requires RQMC points in larger dimension. We examine in particular the variance-gamma, Heston, and Ornstein-Uhlenbeck stochastic volatility models, and we provide numerical results.

\section{INTRODUCTION}
\label{sec:intro}

Quasi-Monte Carlo (QMC) and randomized QMC (RQMC) methods can improve efficiency significantly 
when estimating an integral in a moderate number of dimensions, but their use for simulating 
Markov chains over a large number of steps has been limited so far.
The array-RQMC method, developed for that purpose, has been show to work well for some 
chains having a low-dimensional state.  
It simulates an array of $n$ copies of the Markov chain so that each chain follows its exact distribution,
but the copies are not independent, and the empirical distribution of the states at any given step
of the chain is a ``low-discrepancy'' approximation of the exact distribution.
At each step, the $n$ chains (or states) are matched one-to-one to a set of $n$ RQMC points
whose dimension is the dimension of the state plus the number of uniform random numbers required
to advance the chain by one more step.  The first coordinates of the points are used to match the 
states to the points and the other coordinates provide the random numbers needed to determine the next state.
When the chains have a large-dimensional state, the dimension used for the match can be reduced via
a mapping to a lower-dimensional space.  Then the matching is performed by sorting both the points and
the chains. When the dimension of the state exceeds 1, this matching is done via a multivariate sort. 
The main idea is to evolve the array of chains in a way that from step to step, the empirical distribution
of the states keeps its low discrepancy.
For further details on the methodology, sorting strategies, convergence analysis, applications, 
and empirical results, we refer the reader to 
\shortciteN{vLEC04m,vDEM05a,vLEC06a,vLEC07b,vLEC08a,vELH08a,vLEC09d,vELH10a,vDIO10a,vLEC10c,tGER15a,vLEC18b}, 
and the other references given there.

The aim of this paper is to examine how Array-RQMC can be applied for option pricing under a
stochastic volatility process such as the variance gamma, Heston, and Ornstein-Uhlenbeck models.
We explain and compare various implementation alternatives, and report empirical experiments to
assess the (possible) gain in efficiency and convergence rate.
A second objective is for the WSC community to become better aware of this method, 
which can have numerous other applications.

Array-RQMC has already been applied for pricing Asian options when the underlying process 
evolves as a geometric Brownian motion (GBM) with fixed volatility \shortcite{vLEC09d,vLEC18b}.
In that case, the state is two-dimensional (it contains the current value of the GBM and its running average)
and a single random number is needed at each step, so the required RQMC points are three-dimensional.
In their experiments, \shortciteN{vLEC18b} observed an empirical variance of the average payoff 
that decreased approximately as $\cO(n^{-2})$ for Array-RQMC, in a range of reasonable values of $n$,
compared with $\cO(n^{-1})$ for independent random points (Monte Carlo).
For $n=2^{20}$ (about one million chains), the variance ratio between Monte Carlo 
and Array-RQMC was around 2 to 4 millions. 

In view of this spectacular success, one wonders how well the method would
perform when the underlying process is more involved, e.g., when it has stochastic volatility.
This is relevant because stochastic volatility models are  more
realistic than the plain GBM model \shortcite{fMAD90a,fMAD98a}.
Success is not guaranteed
because the dimension of the required RQMC points is larger.
For the Heston model, for example, the RQMC points must be five-dimensional instead of three-dimensional,
because the state has three dimensions and we need two uniform random numbers at each step.
It is unclear a priori if there will be any significant variance reduction for reasonable values of $n$.

The remainder is organized as follows. 
In Section~\ref{sec:background}, we state our general Markov chain model and provide background 
on the Array-RQMC algorithm, including matching and sorting strategies.
In Section \ref{sec:experimental}, we describe our experimental setting, 
and the types of RQMC point sets that we consider.
Then we study the application of Array-RQMC under the variance-gamma model in Section~\ref{sec:vg},
the Heston model in Section~\ref{sec:heston}, and
the Ornstein-Uhlenbeck model in Section~\ref{sec:ou}.
We end with a conclusion.


\section{BACKGROUND: MARKOV CHAIN MODEL, RQMC, AND ARRAY-RQMC} 
\label{sec:background}

The option pricing models considered in this paper fit the following framework,
which we use to summarize the Array-RQMC algorithm.
We have a  discrete-time \emph{Markov chain} $\{X_j,\, j\ge 0\}$ defined by a stochastic recurrence	over a measurable 
state space $\cX$:
\begin{equation}
\label{eq:markov-chain}
	X_{0} = x_{0},   \qquad\text{and}\qquad  X_{j} = \varphi_{j} (X_{j-1},\bU_{j}), \quad j = 1,\dots,\tau. 
\end{equation}
where $x_0\in\cX$ is a deterministic initial state,
$\bU_{1},\bU_{2},...$ are independent random vectors uniformly distributed over the
$d$-dimensional unit cube $(0,1)^{d}$,     
the functions $\varphi_{j} : \cX \times(0,1)^{d}\rightarrow\cX$ are measurable, 
and $\tau$ is a fixed positive integer (the time horizon). 
The goal is:
\[
  \mbox{Estimate} \quad \mu_{\rmy} = \EE[Y], \qquad\mbox{ where \ } Y = g(X_{\tau}) 
\]
and $g : \cX \to \RR$ is a \emph{cost} (or reward) function.
Here we have a cost only at the last step 
but in general there can be a cost function for each step and $Y$ would be the sum of 
these costs \shortcite{vLEC08a}.

\emph{Crude Monte Carlo} estimates $\mu$ by the average
$
  \bar Y_n = \frac{1}{n} \sum_{i=0}^{n-1} Y_{i},
$
where $Y_0,\dots,Y_{n-1}$ are $n$ independent realizations of $Y$.
One has $\EE[{\bar Y_n}] =  \mu_{\rmy}$ and  
$\Var[{\bar Y_n}] = \Var[Y]/n$, assuming that $\EE[Y^2] = \sigma^2_{\rmy} < \infty$.
Note that the simulation of each realization of $Y$ requires a vector 
$\bV = (\bU_1,\dots,\bU_\tau)$ of $d\tau$ independent uniform random variables over $(0,1)$,
and crude Monte Carlo produces $n$ independent replicates of this random vector.

\emph{Randomized quasi-Monte Carlo} (RQMC) replaces the $n$ independent realizations 
of $\bV$ by $n$ \emph{dependent} realizations, which form an RQMC point set in $d\tau$ dimensions.
That is, each $\bV_i$ has the uniform distribution over $[0,1)^{d\tau}$, and 
the point set $P_{n}=\{V_{0},...,V_{n-1}\} $ covers $[0,1)^{d\tau}$ more evenly than typical 
independent random points.
With RQMC, ${\bar Y_n}$ remains an unbiased estimator of $\mu$, but its variance can be much smaller,
and can converge faster than $\cO(1/n)$ under certain conditions.
For more details, see \shortciteN{rDIC10a,vLEC00b,vLEC09f,vLEC18a}, for example.
However, when $d\tau$ is large, standard RQMC typically becomes ineffective, in the sense that it
does not bring much variance reduction unless the problem has special structure.
	
\emph{Array-RQMC} is an alternative approach developed specifically for Markov chains
\shortcite{vLEC06a,vLEC08a,vLEC18b}.
To explain how it works, let us first suppose for simplicity (we will relax it later) that there is a mapping 
$h : \cX \to\RR$, that assigns to each state a \emph{value} (or score) which summarizes in a single 
real number the most important information that we should retain from that state
(like the value function in stochastic dynamic programming).  
This $h$ is called the \emph{sorting function}. 
The algorithm simulates $n$ (dependent) realizations of the chain ``in parallel''.
Let $X_{i,j}$ denote the state of chain $i$ at step $j$, for $i=0,\dots,n-1$ and $j=0,\dots,\tau$.
At step $j$, the $n$ chains are sorted by increasing order of their values of $h(X_{i,j-1})$,
the $n$ points of an RQMC point set in $d+1$ dimensions are sorted by their first coordinate,
and each point is matched to the chain having the same position in this ordering.
Each chain $i$ is then moved forward by one step, from state $X_{i,j-1}$ to state $X_{i,j}$, 
using the $d$ other coordinates of its assigned RQMC point.  
Then we move on to the next step, the chains are sorted again, and so on.

The sorting function can in fact be more general and have the form $h : \cX \to\RR^c$
for some small integer $c \ge 1$. Then the mapping between the chains and the points must be
realized in a $c$-dimensional space, i.e., via some kind of $c$-dimensional multivariate sort. 
The RQMC points then have $c+d$ coordinates, and are sorted with the same 
$c$-dimensional multivariate sort based on their first $c$ coordinates, and mapped to the 
corresponding chains.  The other $d$ coordinates are used to move the chains ahead by one step.
In practice, the first $c$ coordinates of the RQMC points do not have to be randomized at each step;
they are usually fixed and the points are already sorted in the correct order based on these coordinates.

Some multivariate sorts are described and compared by \shortciteN{vELH08a,vLEC09d,vLEC18a}.
For example, in a \emph{multivariate batch sort}, we select positive integers 
$n_1, \dots, n_c$ such that $n = n_1 \dots n_c$. 
The states are first sorted by their first coordinate in ${n_1}$ packets of size ${n/n_1}$,
then each packet is sorted by the second coordinate into ${n_2}$ packets of size ${n/n_1 n_2}$, and so on.
The RQMC points are sorted in exactly the same way, based on their first $c$ coordinates.
In the \emph{multivariate split sort}, we assume that $n = 2^e$
and we take $n_1 = n_2 = \cdots = n_e = 2$.  
That is, we first split the points in 2 packets based on the first coordinate,
then split each packet in two by the second coordinate, and so on.
If $e > c$, after $c$ splits we get back to the first coordinate and continue.

Examples of heuristic sorting functions $h : \cX\to\RR$ are given in \shortcite{vLEC08a,vLEC18b}.
\citeN{vWAC08a} and \citeN{tGER15a} suggested to first map the $c$-dimensional states to $[0,1]^c$
and then use a space filling curve in $[0,1]^c$ to map them to $[0,1]$, which provides a total order.
\citeN{tGER15a} proposed to map the states to $[0,1]^c$ via a component-wise rescaled logistic transformation,
then order them with a Hilbert space-filling curve.   See \shortciteN{vLEC18b} for a more detailed discussion.
Under smoothness conditions, they proved that the resulting unbiased Array-RQMC estimator has $o(1/n)$ variance,
which beats the $\cO(1/n)$ Monte Carlo rate. 

Algorithm~\ref{algo:array-rqmc} states the Array-RQMC procedure in our setting. 
Indentation delimits the scope of the \textbf{for} loops.
For any choice of sorting function $h$, the average $\hat\mu_{\arqmc,n} = \bar Y_n$ 
returned by this algorithm is always an \emph{unbiased} estimator of $\mu$.
An unbiased estimator of $\Var[\bar Y_n]$ can be obtained by making $m$ independent realizations 
of $\hat\mu_{\arqmc,n}$ and computing their \emph{empirical variance}.

\begin{algorithm}[ht]
\caption{: Array-RQMC Algorithm for Our Setting}
\label{algo:array-rqmc}
\begin{algorithmic}
\State {\textbf{for} $i=0,\dots,n-1$ \textbf{do} ${X_{i,0}} \g x_0$;}
\For {${j} = 1, 2, \dots, \tau$}
  \State Sorting: Compute an appropriate permutation $\pi_{j}$ of the $n$ chains, based on 
	\State $\qquad$ the $h(X_{i,j-1})$, to match the $n$ states with the RQMC points;
  \State Randomize afresh 
      the RQMC points $\{\bU_{0,j},\dots,\bU_{n-1,j}\}$;
  \State {\textbf{for} $i=0,\dots,n-1$ \textbf{do} ${X_{i,j}} = \varphi_j(X_{\pi_j(i),j-1}, \bU_{i,j})$;}
\EndFor
\Return the average $\hat\mu_{\arqmc,n} = \bar Y_n = (1/n) \sum_{i=0}^{n-1} g(X_{i,\tau})$ 
  as an estimate of $\mu_{\rmy}$.
\end{algorithmic}
\end{algorithm}

\section{EXPERIMENTAL SETTING}
\label{sec:experimental}

For all the option pricing examples in this paper, we have an asset price that evolves 
as a stochastic process $\{S(t),\, t\ge 0\}$ and a payoff that depends on the values of 
this process at fixed observation times $0 = t_{0} <t_{1}<t_{2}< ... < t_{c}=T$.
More specifically, for given constants $r$ (the interest rate) and $K$ (the strike price),  
we consider an \emph{European option} whose payoff is 
\[
  Y = Y_{\rme} = g(S(T)) = e^{-rT} \max (S(T)-K, 0)
\]
and a discretely-observed \emph{Asian option} whose payoff is 
\[
  Y = Y_{\rma} = g(\bar S) = e^{-rT} \max (\bar S-K, 0)
\]
where $\bar S = (1/c) \sum _{j=1}^{c} S(t_{j})$.
In this second case, the running average  
$\bar S_j = (1/j)\sum_{\ell=1}^j S(t_{\ell})$ must be kept in the state of the Markov chain.
The information required for the evolution of $S(t)$ depends on the model and is given 
for each model in forthcoming sections.  It must be maintained in the state.  
For the case where $S$ is a plain GBM, the state of the Markov chain
at step $j$ can be taken as $X_j = (S(t_j), \bar S_j)$, a two-dimensional state, as was done in 
\shortciteN{vLEC09d} and \shortciteN{vLEC18b}.

In our examples, the states are always multidimensional.  
To match them with the RQMC points, we will use a split sort, a batch sort, and a Hilbert-curve sort,
and compare these alternatives.  The Hilbert sort requires a transformation of the $\ell$-dimensional states
to the unit hypercube $[0,1]^\ell$.  For this, we use a \emph{logistic transformation} defined by
$\psi(x)=(\psi_{1}(x_{1}),...,\psi_{\ell}(x_{\ell})) \in[0,1]^\ell$ for all $x = (x_1,\dots,x_\ell) \in\cX$, where
\begin{equation} 
	\psi_{j}(x_{j})= 
	  \left [1+\exp \left(-\frac{x_{j}-\underline{x}_{j}}{\bar{x}_{j}-\underline{x}_{j}} \right) \right]^{-1}, 
		\quad  j=1,...,\ell,
\end{equation}
with constants $\bar{x}_{j}=\mu_{j}+2 \sigma_{j} $ and $\underline{x}_{j}=\mu_{j}-2 \sigma_{j}$
in which $ \mu_{j} $ and $ \sigma_{j}$ are estimates of the mean and the variance of the distribution 
of the $j$th coordinate of the state.
In Section~\ref{sec:vg}, we will also consider just taking a linear combination 
of the two coordinates, to map a two-dimensional state to one dimension.

For RQMC, we consider 
\vskip-25pt\null%
\begin{verse}
(1) Independent points, which corresponds to crude Monte Carlo (MC);\\
(2) Stratified sampling over the unit hypercube (Stratif);\\
(3) Sobol' points with a random linear matrix scrambling 
    and a digital random shift (Sobol+LMS);\\
(4) Sobol' points with nested uniform scrambling (Sobol+NUS);\\
(5) A rank-1 lattice rule with a random shift modulo 1 followed by a baker's transformation (Lattice+baker).
\end{verse}
\vskip-12pt
The first two are not really RQMC points, but we use them for comparison.
For stratified sampling, we divide the unit hypercube into $n = k^{\ell+d}$ congruent subcubes
for some integer $k > 1$, and we draw one point randomly in each subcube.
For a given target $n$, we take $k$ as the integer for which $k^{\ell+d}$ is closest to this target $n$.
For the Sobol' points, we took the default direction numbers in SSJ, 
which are from \shortciteN{iLEM04a}.
The LMS and NUS randomizations are explained in \citeN{vOWE03a} and \citeN{vLEC09f}.
For the rank-1 lattice rules, we used generating vectors found by Lattice Builder \shortcite{vLEC16a},
using the $\cP_2$ criterion with order-dependent weights $(0.8)^k$ for projections of order $k$.

For each example, each sorting method, each type of point set, and each selected value of $n$,
we ran simulations to estimate $\Var[\bar Y_n]$.   
For the stratified and RQMC points, this variance was estimated by replicating the RQMC scheme $m = 100$ times
independently.  
For a fair comparison with the MC variance $\sigma^2_{\rmy} = \Var[Y]$, for these point sets we used the 
\emph{variance per run}, defined as $n \Var[\bar Y_n]$.
We define the \emph{variance reduction factor} (VRF) for a given method compared with MC 
by $\sigma^2_{\rmy} / (n \Var[\bar Y_n])$.
In each case, we fitted a linear regression model for the variance per run as a function of $n$,
in log-log scale.  We denote by $\hat\beta$ the regression slope estimated by this linear model.

In the remaining sections, we explain how the process $\{S(t),\, t\ge 0\}$ is defined in each
case, how it is simulated.  We show how we can apply Array-RQMC and we provide numerical results.
All the experiments were done in Java using the SSJ library \cite{sLEC05a,iLEC16j}.

\section{OPTION PRICING UNDER A VARIANCE-GAMMA PROCESS}
\label{sec:vg}

The \emph{variance-gamma} (VG) model was proposed for option pricing by \citeN{fMAD90a} and \citeN{fMAD98a},
and further studied by \shortciteN{fFU98a,fAVR03a,fAVR06a}, for example.
A VG process is essentially a Brownian process for which the time clock runs
at random and time-varying speed driven by a gamma process.
The VG process with parameters $(\theta,\sigma^2,\nu)$ 
is defined as $Y = \{Y(t) = X(G(t)),\, t\ge 0\}$ where 
$X = \{X(t),\, t\ge 0\}$ is a Brownian motion with drift and variance parameters $\theta$ and $\sigma^2$, and
$G = \{G(t),\, t\ge 0\}$ is a gamma process with drift and volatility parameters $1$ and $\nu$,
independent of $X$.   
This means that $X(0) = 0$, $G(0)=0$, both $B$ and $G$ have independent increments,
and for all $ t \geq 0 $ and  $ \delta > 0 $, we have 
$X(t+\delta)-X(t)\sim {\rm Normal}(\delta\theta,\delta\sigma^{2})$,
a normal random variable with mean $\delta\theta$ and variance $\delta\sigma^{2}$,
and $G(t+\delta)-G(t) \sim {\rm Gamma}(\delta/\nu, \nu)$, a gamma random variable
with mean $\delta$ and variance $\delta \nu$.
The gamma process is always non-decreasing, which ensures that the time clock never goes backward.
In the VG model for option pricing, the asset value follows the \emph{geometric variance-gamma} (GVG) 
process  $S = \{S(t),\, t\ge 0\}$ defined by
\[
  S(t) = S(0) \exp\left[ (r+\omega) t + X(G(t)) \right],
\]
where $\omega = \ln (1-\theta\nu-\sigma^2 \nu/2)/\nu$.

To generate realizations of $\bar S$  
for this process, we must generate $S(t_1),\dots,S(t_{\tau})$, and there are many ways of doing this. 
With Array-RQMC, we want to do it via a Markov chain with a low-dimensional state.
The running average $\bar S_j$ must be part of the state, as well as sufficient information
to generate the future of the path.
A simple procedure for generating the path is to sample sequentially $G(t_1)$,
then $Y(t_1) = X(G(t_1))$ conditional on $G(t_1)$, 
then $G(t_2)$ conditional on $G(t_1)$, 
then $Y(t_2) = X(G(t_2))$ conditional on $(G(t_1), G(t_2), Y(t_1))$, and so on.
We can then compute any $S(t_j)$ directly from $Y(t_j)$.

It is convenient to view the sampling of $(G(t_j),\, Y(t_j))$ conditional on $(G(t_{j-1}),\, Y(t_{j-1}))$
as one step (step $j$) of the Markov chain.
The state of the chain at step $j-1$ can be taken as $X_{j-1} = (G(t_{j-1}),\, Y(t_{j-1}),\, \bar S_{j-1})$,
so we have a three-dimensional state, and we need two independent uniform random numbers at each step,
one to generate $G(t_j)$ and the other to generate $Y(t_j) = X(G(t_j))$ given $(G(t_{j-1}), G(t_j), Y(t_{j-1}))$,
both by inversion. 
Applying Array-RQMC with this setting would require a five-dimensional RQMC point set at each step,
unless we can map the state to a lower-dimensional representation.

However, a key observation here is that the distribution of the increment 
$\Delta Y_j = Y(t_j) - Y(t_{j-1})$ depends only on the increment $\Delta_j = G(t_j) - G(t_{j-1})$
and not on $G(t_{j-1})$.  This means that there is no need to memorize the latter in the state!
Thus, we can define the state at step $j$ as the two-dimensional vector $X_{j} = (Y(t_{j}),\, \bar S_{j})$,
or equivalently $X_{j} = (S(t_{j}),\, \bar S_{j})$,
and apply Array-RQMC with a four-dimensional RQMC point set if we use a two-dimensional sort for the states,
and a three-dimensional RQMC point set if we map the states to a one-dimensional representation
(using a Hilbert curve or a linear combination of the coordinates, for example).
At step $j$, we generate $\Delta_{j} \sim \dGamma((t_{j}-t_{j-1})/\nu, \nu)$ by inversion using a 
uniform random variate $U_{j,1}$, i.e., via $\Delta_{j} = F_j^{-1}(U_{ j,1})$ where $F_j$ is the cdf
of the $\dGamma((t_{j}-t_{j-1})/\nu, \nu)$ distribution, then $\Delta Y_j$ by inversion from the 
normal distribution with mean $\theta \Delta_{j}$ and variance $\sigma^2 \Delta_{j}$, 
using a uniform random variate $U_{j,2}$.
Algorithm~\ref{algo:sequential} summarizes this procedure.
The symbol $\Phi$ denotes the standard normal cdf.
We have 
\[
  X_j = (Y(t_{j}),\bar{S}_{j}) = \varphi_{j}(Y(t_{j-1}), \bar{S}_{j-1}, U_{j,1}, U_{j,2}) 
\]
where $\varphi_j$ is defined by the algorithm.
The payoff function is $g(X_c) = \bar S_c = \bar S$.

\begin{algorithm}
\caption{Computing $X_j = (Y(t_{j}),\bar{S}_{j})$ given $(Y(t_{j-1}), \bar{S}_{j-1})$, for $1\leq j\leq \tau$.} 
\label{algo:sequential}
\begin{algorithmic}
			\State {Generate \ $U_{j,1}, U_{j,2} \sim \dUnif(0,1)$, independent;}
			\State {$\Delta_{j} = F_j^{-1}(U_{j,1}) \sim {\rm Gamma}((t_{j}-t_{j-1})/\nu, \nu)$;}
			\State {$Z_{j} = \Phi^{-1}(U_{j,2}) \sim \dNormal(0,1)$;}
			\State {$Y(t_{j})\leftarrow  Y(t_{j-1}) + \theta\Delta_{j}+ \sigma \sqrt{\Delta_{j}} Z_j$;}
			\State {$S(t_j) \leftarrow  S(0) \exp[(r+\omega) t_j + Y(t_j)]$;}
			\State {$\bar S_j = [(j-1) \bar S_{j-1} + S(t_j)]/j$;}
\end{algorithmic}
\end{algorithm}

With this two-dimensional state representation, if we use a split sort or batch,
we need four-dimensional RQMC points.
With the Hilbert-curve sort, we only need three-dimensional RQMC points.
We also tried a simple linear mapping $h_j : \RR^2\to\RR$ defined by 
$h_j(S(t_j), \bar S_j) = b_j \bar S_j + (1-b_j) S(t_j)$
where $b_j = (j-1)/(\tau-1)$.
At each step $j$, this $h_j$ maps the state $X_j$ to a real number $h_j(X_j)$,
and we sort the states by increasing order of their value of $h_j(X_j)$.
It uses a convex linear combination of $S(t_j)$ and $\bar S_j$ whose coefficients depend on $j$.
The rationale for the (heuristic) choice of $b_j$  
is that in the late steps (when $j$ is near $\tau$), the current average $\bar{S}_{j}$
is more important (has more predictive power for the final payoff) than the current $S(t_{j})$,
whereas in the early steps, the opposite is true.  

We made an experiment with the following model parameters, taken from \citeN{fAVR06a}: 
$\theta = -0.1436 $, $ \sigma = 0.12136 $, $\nu = 0.3 $, $r = 0.1 $, $T = 240/365$, $\tau = 10$, 
$t_j = 24j/365$ for $j=1,\dots,\tau$,  $K = 100$, and $S(0) = 100$.   
The time unit is one year, the horizon is 240 days, and there is an observation time every 24 days.
The exact value of the expected payoff for the Asian option is $\mu \approx 8.36 $,
and the MC variance per run is $\sigma_{\rmy}^{2} = \Var[Y_{\rma}] \approx 59.40 $.

\begin{table}[ht ]
		\centering
		\small
\caption{Regression slopes $\hat\beta$ for $\log_{2} \Var[\mua]$ vs $\log_2(n)$, 
    and VRF compared with MC for $n=2^{20}$, denoted VRF20, for the Asian option under the VG model }
\label{tab:vg-array} 
\begin{tabular}{l|r|c|r}
			\hline
		Sort & Point sets & $ \hat{\beta} $ & VRF20  \\
		\hline
		\multirow{4}{2cm}{Split sort }	& MC & -1 &  1  \\
		& Stratif & -1.17 &  42  \\
		& Sobol'+LMS & -1.77  &  91550  \\
		& Sobol'+NUS & -1.80  &   106965   \\
		& Lattice+baker & -1.83 &    32812  \\
		\hline 	     	    
		\multirow{4}{2cm}{Batch sort ($ n_{1}=n_{2} $) }	 	& MC & -1 &  1  \\
		& Stratif & -1 &  42 \\
		& Sobol'+LMS & -1.71 &   100104   \\
		& Sobol'+NUS & -1.54 &  90168    \\
		& Lattice+baker & -1.95  &  58737    \\
		\hline 		
		\multirow{4}{2cm}{Hilbert sort (with logistic map) }	 	& MC & -1 &  1  \\
		& Stratif & -1.43 & 204  \\
		& Sobol'+LMS & -1.59 &  68297   \\
		& Sobol'+NUS & -1.67 &    79869   \\
		& Lattice+baker & -1.55  & 45854   \\
		\hline 	
		\multirow{4}{2cm}{Linear map sort  }	 	& MC & -1 & 1   \\
		& Stratif & -1.35 & 192  \\
		& Sobol'+LMS & -1.64 &   115216  \\
		& Sobol'+NUS & -1.75 &   166541    \\
		& Lattice+baker & -1.72  &  68739   \\
		\hline 			
	\end{tabular} 	
\end{table}

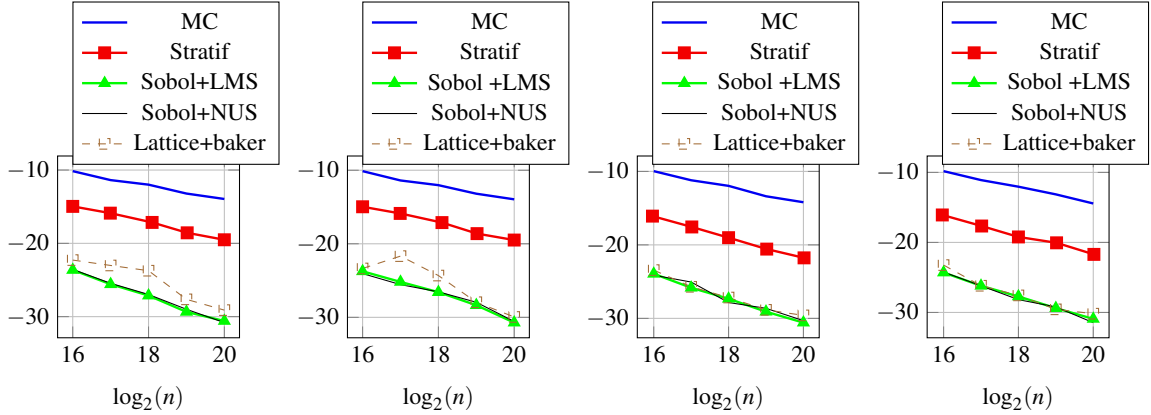
\begin{figure}[!hbt] \footnotesize  
\centering
\begin{tikzpicture} 
	\begin{axis}[ 
		cycle list name=defaultcolorlist,
		legend style={at={(1.22,1.85)}, anchor={north east}},
		xlabel=$\log_2(n)$,
		grid,
		width=4cm,
		height=4cm,
		] 
		\addplot table[x=log(n),y=log(MISE)] { 
		log(n)  log(MISE) 
		16.0  -10.151201808434037 
		17.0  -11.374817074545694 
		18.0  -11.994381098726798 
		19.0  -13.208370282151458 
		20.0  -13.949334962385464 
	}; 
	\addlegendentry{MC}
	\addplot table[x=log(n),y=log(MISE)] { 
		log(n)  log(MISE) 
		16.0  -14.959378903560621 
		16.99171005377434  -15.88089020607308 
		18.094247824228052  -17.14658716084826 
		19.019550008653873  -18.564637351701542 
		20.0  -19.513185167970356 
	}; 
	\addlegendentry{Stratif}
	
	\addplot table[x=log(n),y=log(MISE)] { 
		log(n)  log(MISE) 
		16.0  -23.633704459633822 
		17.0  -25.580879149864206 
		18.0  -27.124986360860927 
		19.0  -29.348723118567527 
		20.0  -30.58977405394141 
	}; 
	\addlegendentry{Sobol+LMS}
	%
	%
	\addplot table[x=log(n),y=log(MISE)] { 
		log(n)  log(MISE) 
		16.0  -23.55879443464881 
		17.0  -25.452065755021675 
		18.0  -26.978058126434465 
		19.0  -28.994972675511697 
		20.0  -30.81428988989172 
	}; 
	\addlegendentry{Sobol+NUS}
	%
	\addplot table[x=log(n),y=log(MISE)] { 
		log(n)  log(MISE) 
		16.0  -22.28254034764946 
		17.0  -23.01761769453436 
		18.0  -23.731399349454925 
		19.0  -27.640591609428213 
		20.0  -29.109451903531756 
		
	}; 
	\addlegendentry{Lattice+baker}
	\end{axis}
		\end{tikzpicture}
		\begin{tikzpicture}  
		\begin{axis}[ 
		cycle list name=defaultcolorlist,
		legend style={at={(1.22,1.85)}, anchor={north east}},
		xlabel=$\log_2(n)$,
		grid,
		width=4cm,
		height=4cm,
		] 
		\addplot table[x=log(n),y=log(MISE)] { 
		log(n)  log(MISE) 
		16.0  -10.125167088824263 
		17.0  -11.391734398723784 
		18.0  -12.030243862257686 
		19.0  -13.19124273647122 
		20.0  -13.957222521780116 
	}; 
	\addlegendentry{MC}

	\addplot table[x=log(n),y=log(MISE)] { 
		log(n)  log(MISE) 
		16.0  -14.967680095056181 
		16.99171005377434  -15.872803090894319 
		18.094247824228052  -17.11791662799947 
		19.019550008653873  -18.618933749651113 
		20.0  -19.491670326261605 
	}; 
	\addlegendentry{Stratif}
	\addplot table[x=log(n),y=log(MISE)] { 
		log(n)  log(MISE) 
		16.0  -23.742038334288367 
		17.0  -25.18529149145997 
		18.0  -26.575981540352245 
		19.0  -28.349253143579467 
		20.0  -30.71863728507717 
	}; 
	\addlegendentry{Sobol +LMS}
	
	%
	%
	\addplot table[x=log(n),y=log(MISE)] { 
		log(n)  log(MISE) 
		16.0  -24.05405256317996 
		17.0  -25.551335850072785 
		18.0  -26.532683096004593 
		19.0  -27.949058481148693 
		20.0  -30.56782525433717 
	}; 
	\addlegendentry{Sobol+NUS}
	%
	\addplot table[x=log(n),y=log(MISE)] { 
		log(n)  log(MISE) 
		16.0  -23.354622638739173 
		17.0  -21.675118811834484 
		18.0  -24.319492844187447 
		19.0  -27.946213398912015 
		20.0  -29.94949071310452 
		
	}; 
	\addlegendentry{Lattice+baker}
	\end{axis}
		\end{tikzpicture}
		\begin{tikzpicture} 
		\begin{axis}[ 
		cycle list name=defaultcolorlist,
		legend style={at={(1.22,1.85)}, anchor={north east}},
		xlabel=$\log_2(n)$,
		grid,
		width=4cm,
		height=4cm,
		] 
		\addplot table[x=log(n),y=log(MISE)] { 
		log(n)  log(MISE) 
		16.0  -9.965506807807301 
		17.0  -11.205344676647975 
		18.0  -11.967883997584597 
		19.0  -13.384159678031219 
		20.0  -14.196724538935994 
	}; 
	\addlegendentry{MC}

	\addplot table[x=log(n),y=log(MISE)] { 
		log(n)  log(MISE) 
		15.965784284662087  -16.08499538939687 
		17.017276025914487  -17.56037697635942 
		18.0  -19.018210419986847 
		19.019550008653873  -20.572247692334237 
		20.017276025914487  -21.7649559791986 
	}; 
	\addlegendentry{Stratif}
	\addplot table[x=log(n),y=log(MISE)] { 
		log(n)  log(MISE) 
		16.0  -23.93571019224644 
		17.0  -25.832152908036722 
		18.0  -27.33008775845543 
		19.0  -29.091177852445593 
		20.0  -30.578431433321636 
	}; 
	\addlegendentry{Sobol +LMS}
	
	%
	\addplot table[x=log(n),y=log(MISE)] { 
		log(n)  log(MISE) 
		16.0  -24.120858868108805 
		17.0  -25.058596678607625 
		18.0  -27.841113569201394 
		19.0  -28.632361952794387 
		20.0  -30.308379158193635 
	}; 
	\addlegendentry{Sobol+NUS}
	%
	%
	\addplot table[x=log(n),y=log(MISE)] { 
		log(n)  log(MISE) 
		16.0  -23.42830483146864 
		17.0  -25.702792296644837 
		18.0  -27.16182898294017 
		19.0  -28.83775332640567 
		20.0  -29.59227100834182 
		
	}; 
	\addlegendentry{Lattice+baker}
	\end{axis}
		\end{tikzpicture}
		\begin{tikzpicture} 
		\begin{axis}[ 
		cycle list name=defaultcolorlist,
		legend style={at={(1.22,1.85)}, anchor={north east}},
		xlabel=$\log_2(n)$,
		grid,
		width=4cm,
		height=4cm,
		] 
		\addplot table[x=log(n),y=log(MISE)] { 
		log(n)  log(MISE) 
		16.0  -9.82305725456555 
      	17.0  -11.099014117169833 
      	18.0  -12.044117067884285 
      	19.0  -13.140553493058876 
      	20.0  -14.440549409485193 
	}; 
	\addlegendentry{MC}

	\addplot table[x=log(n),y=log(MISE)] { 
		log(n)  log(MISE) 
      	15.965784284662087  -16.08023037439075 
      	17.017276025914487  -17.658391548343992 
      	18.0  -19.220332169166976 
      	19.019550008653873  -20.058388378525507 
      	20.017276025914487  -21.710452609694315 
	}; 
	\addlegendentry{Stratif}
	\addplot table[x=log(n),y=log(MISE)] { 
		log(n)  log(MISE) 
		16.0  -24.301030870661574 
      	17.0  -26.19935628004586 
		18.0  -27.716727137785195 
		19.0  -29.385458064674975 
		20.0  -30.921744935733827 
	}; 
	\addlegendentry{Sobol +LMS}
	
	\addplot table[x=log(n),y=log(MISE)] { 
		log(n)  log(MISE) 
      	16.0  -24.23697989879872 
      	17.0  -26.19235071537741 
      	18.0  -28.169357039010745 
      	19.0  -29.216652851994684 
      	20.0  -31.453280459358083 
	}; 
	\addlegendentry{Sobol+NUS}
	%
	%
	\addplot table[x=log(n),y=log(MISE)] { 
		log(n)  log(MISE) 
      	16.0  -23.236130419924017 
      	17.0  -26.186668822233617 
      	18.0  -27.56543498454785 
      	19.0  -29.55463426017593 
      	20.0  -30.17659115439563
	}; 
	\addlegendentry{Lattice+baker}
	\end{axis}
\end{tikzpicture}
\caption{Plots of empirical $\log_{2}\Var[\mua]$ vs $\log_2(n)$ for various sorts and point sets, 
  based on $m=100$ independent replications.
	Left to right: split sort, batch sort, Hilbert sort, linear map sort.}
\label{fig:vg-array}
\end{figure}

Table~\ref{tab:vg-array} summarizes the results.
For each selected sorting method and point set, we report the estimated slope $\hat\beta$ for the 
linear regression model of $\log_{2}\Var[\mua]$ as a function of $\log_2(n)$ 
obtained from $m=100$ independent replications with $n=2^{e}$ for $e=16,...,20$,
as well as the variance reduction factors (VRF) observed for $n=2^{20}$ (about one million samples), 
denoted VRF20.
For MC, the exact slope (or convergence rate) $\beta$ is known to be $\beta = -1$. 
We see from the table that Array-RQMC provides much better convergence rates (at least empirically),
and reduces the variance by very large factors for $n = 2^{20}$.
Interestingly, the largest factors are obtained with the Sobol' points combined with our 
heuristic linear map sort, although the other sorts are also doing quite well.
Figure~\ref{fig:vg-array} shows plots of $\log_{2}\Var[\mua]$ vs $\log_2(n)$ for selected sorts.
It gives an idea of how well the linear model fits in each case.

There are other ways of defining the steps of the Markov chain for this example.
For example, one can have one step for each $\dUnif(0,1)$ random number 
that is generated.  This would double the number of steps, from $c$ to $2c$.
We generate $\Delta_1$ in the first step, $Y(t_1)$ in the second step, 
$\Delta_2$ in the third step, $Y(t_2)$ in the fourth step, and so on.
Generating a single uniform per step instead of two reduces by 1 the dimension of
the required RQMC point set.  At odd step numbers, when we generate a $\Delta_j$,
the state can still be taken as $(Y(t_{j-1}), \bar{S}_{j-1})$ and we only need three-dimensional
RQMC points, so we save one dimension.
But at even step numbers, we need $\Delta_j$ to generate $Y(t_j)$,
so we need a three-dimensional state  $(Y(t_{j-1}), \Delta_j, \bar{S}_{j-1})$
and four-dimensional RQMC points.
We tried this approach and it did not perform better than the one described earlier,
with two uniforms per step.  It is also more complicated to implement.

\shortciteN{fAVR03a,fAVR06a} describe other ways of simulating the VG process, for instance
Brownian and gamma bridge sampling (BGBS) and difference of gammas bridge sampling (DGBS).
BGBS generates first $G(t_c)$ then $Y(t_c)$, then conditional on this it generates 
$G(t_{c/2})$ then $Y(t_{c/2})$ (assuming that $c$ is even), and so on.
DGBS writes the VG process $Y$ as a difference of two independent gamma processes 
and simulate both using the bridge idea just described: first generate the values of the two 
gamma processes at $t_{c}$, then at $t_{c/2}$, etc.
When using classical RQMC, these sampling methods brings an important variance reduction 
compared with the sequential one we use here for our Markov chain. 
With Array-RQMC, however, they are impractical, because the dimension of the state 
(the number of values that we need to remember) grows up to about $c$, which is much to high,
and the implementation is much more complicated.

For this VG model, we do not report results on the European option with Array-RQMC,
because the Markov chain would have only one step: We can generate directly 
$G(t_c)$ and then $Y(t_c)$. For this, ordinary RQMC works well enough \shortcite{vLEC18a}.

\section{OPTION PRICING UNDER THE HESTON VOLATILITY MODEL}
\label{sec:heston}

The Heston volatility model is defined by the following two-dimensional stochastic differential equation:
	\begin{eqnarray*}
	\d S(t) &=& rS(t)\d t+ V(t)^{1/2}S(t)\d B_{1}(t),\\
	\d V(t) &=& \lambda(\sigma^{2}-V(t))\d t + \xi V(t)^{1/2}\d B_{2}(t),
	\end{eqnarray*} 
for $t \geq 0$, where $(B_{1},B_{2})$ is a pair of standard Brownian motions with correlation 
$\rho $ between  them, $ r $ is the risk-free rate, $ \sigma^{2} $ is the long-term average variance parameter, 
$\lambda $ is the rate of return to the mean for the variance, 
and $\xi$ is a volatility parameter for the variance.
The processes $S = \{S(t),\, t\ge 0\}$ and $V = \{V(t),\, t\ge 0\}$ represent the asset price 
and the volatility, respectively, as a function of time.
We will examine how to estimate the price of European and Asian options with Array-RQMC under this model.
Since we do not know how to generate $(S(t+\delta),\, V(t+\delta))$ exactly from its conditional 
distribution given $(S(t),\, V(t))$ in this case, we have to discretize the time. 
For this, we use the Euler method with $\tau$ time steps of length $\delta = T/\tau$ to generate a skeleton
of the process at times $w_j = j\delta$ for $j=1,\dots,\tau$, over $[0,T]$.
For the Asian option, we assume for simplicity that the observation times $t_1,\dots,t_c$ 
used for the payoff are all multiples of $\delta$, so each of them is equal to some $w_j$.
\hpierre{Here we should avoid the notation $t_j$, because $j$ is the step number of the Markov chain,
  so we would have two different meanings for the index $j$.}

\begin{table}[hbt]
		\centering
			\small
\caption{Regression slopes $\hat\beta$ for $\log_{2} \Var[\mua]$ vs $\log_2(n)$, 
    and VRF compared with MC for $n=2^{20}$, denoted VRF20, for the Asian option under the Heston model. }
\label{tab:heston-array} 
\begin{tabular}{l|r|r|r|r|r}
			\hline
			\multicolumn{2}{c|}{}  & \multicolumn{2}{c|}{European}& \multicolumn{2}{c}{Asian}   \\
			\hline
		Sort &  Point sets & $ \hat{\beta} $ & VRF20 & $ \hat{\beta} $ & VRF20  \\
		\hline
		\multirow{5}{*}{Split sort }	& MC & -1 & 1 &  -1 &  1 \\
		& Stratif & -1.26 & 103  &  -1.29 & 38  \\
		& Sobol'+LMS & -1.59 & 44188   & -1.48 &  6684 \\
		& Sobol'+NUS & -1.46&  30616  & -1.46 &  5755  \\
		& Lattice+baker & -1.50  &  26772 &  -1.55 & 5140   \\   
		\hline 	     	    
		\multirow{5}{2.5cm}{Batch sort 
		}	 	& MC & -1 & 1 &  -1 & 1 \\
		& Stratif & -1.24 & 91 &   -1.25 &  33 \\
		& Sobol'+LMS & -1.66  &  22873   & -1.23 & 815  \\
		& Sobol'+NUS & -1.72  & 30832   & -1.38& 1022  \\
		& Lattice+baker & -1.75 &  12562    & -1.22& 762  \\
		\hline
		\multirow{5}{2.5cm}{Hilbert  sort  (with logistic map)
		}	 	& MC & -1& 1 &  -1 & 1  \\
		& Stratif & -1.26 & 43 &   -1.05 &  29 \\
		& Sobol'+LMS & -1.14  &  368   & -0.87 &  39 \\
		& Sobol'+NUS & -1.06  & 277   & -1.11&  49 \\
		& Lattice+baker & -1.12 &  250    & -0.89 &  42\\
		\hline
	\end{tabular} 
\end{table}

	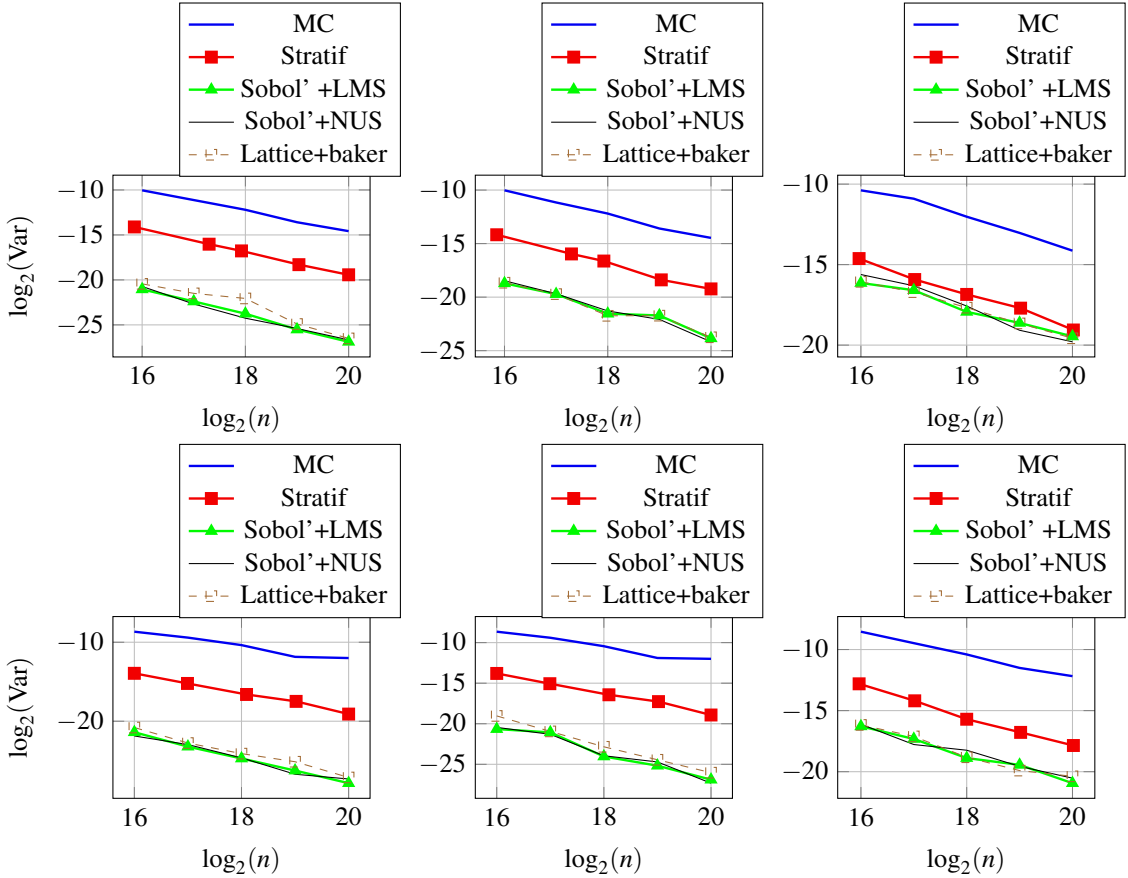
\begin{figure}[htb] \footnotesize 
		\centering
		\small
		\begin{tikzpicture} %
		\begin{axis}[ 
		cycle list name=defaultcolorlist,
		legend style={at={(1.12,1.95)}, anchor={north east}},
		xlabel=$\log_2(n)$,
		ylabel=$\log_2(\Var)$,
		grid,
		width=5cm,
		height=4cm,
		] 
	\addplot table[x=log(n),y=log(MISE)] { 
		log(n)  log(MISE) 
		16.0  -10.0393965896146 
		17.0  -11.122034035538492 
		18.0  -12.203482113338744 
		19.0  -13.588402242719132 
		20.0  -14.578709496795367 
	}; 
	\addlegendentry{MC}

	\addplot table[x=log(n),y=log(MISE)] { 
		log(n)  log(MISE) 
		15.84962500721156  -14.111302800188136 
		17.297158093186486  -16.028321372198608 
		17.92481250360578  -16.774827814961732 
		19.03677461028802  -18.304195272754914 
		20.0  -19.43172105454755 
	}; 
	\addlegendentry{Stratif}
	\addplot table[x=log(n),y=log(MISE)] { 
		log(n)  log(MISE) 
		16.0  -21.050533059699816 
		17.0  -22.395066360734926 
		18.0  -23.75575975421682 
		19.0  -25.512215835716276 
		20.0  -26.881037454834505 
	}; 
	\addlegendentry{Sobol' +LMS}
	%
	\addplot table[x=log(n),y=log(MISE)] { 
		log(n)  log(MISE) 
		16.0  -20.730649960589506 
		17.0  -22.67690957689726 
		18.0  -24.255429521603652 
		19.0  -25.40529102791617 
		20.0  -26.665055006634283 
	}; 
	\addlegendentry{Sobol'+NUS}
	%
	\addplot table[x=log(n),y=log(MISE)] { 
		log(n)  log(MISE) 
		
		16.0  -20.438241859943144 
		17.0  -21.49545519899545 
		18.0  -22.060176245547 
		19.0  -24.903082129918214 
		20.0  -26.502180414676182 
	}; 
	\addlegendentry{Lattice+baker}
	\end{axis}
		\end{tikzpicture}
		\begin{tikzpicture} 
		\begin{axis}[ 
		cycle list name=defaultcolorlist,
		legend style={at={(1.12,1.95)}, anchor={north east}},
		xlabel=$\log_2(n)$,
		grid,
		width=5cm,
		height=4cm,
		] 
		\addplot table[x=log(n),y=log(MISE)] { 
		log(n)  log(MISE) 
		16.0  -10.024845458352644 
		17.0  -11.159811274387627 
		18.0  -12.195612913788786 
		19.0  -13.590093217796193 
		20.0  -14.467680580335442 
		
	}; 
	\addlegendentry{MC}
	\addplot table[x=log(n),y=log(MISE)] { 
		log(n)  log(MISE) 
		15.84962500721156  -14.174715605814937 
		17.297158093186486  -15.965674741102829 
		17.92481250360578  -16.629129192425264 
		19.03677461028802  -18.38704490691805 
		20.0  -19.230981344472582 
		
	}; 
	\addlegendentry{Stratif}
	
	\addplot table[x=log(n),y=log(MISE)] { 
		log(n)  log(MISE) 
		16.0  -18.710312874636735 
		17.0  -19.7195535149438 
		18.0  -21.545707235411527 
		19.0  -21.710845947842163 
		20.0  -23.844994987304258 
		
	}; 
	\addlegendentry{Sobol'+LMS}		
	%
	%
	\addplot table[x=log(n),y=log(MISE)] { 
		log(n)  log(MISE) 
		16.0  -18.471130506776717 
		17.0  -19.687071984460964 
		18.0  -21.2895781658613 
		19.0  -22.077391393023678 
		20.0  -24.171378533597295 
		
	}; 
	\addlegendentry{Sobol'+NUS}
	%
	%
	\addplot table[x=log(n),y=log(MISE)] { 
		log(n)  log(MISE) 
		
		16.0  -18.65403828216375 
		17.0  -19.716095401701306 
		18.0  -21.742953862467 
		19.0  -21.729154800441403 
		20.0  -23.7487392043938 
		
	}; 
	\addlegendentry{Lattice+baker}
	\end{axis}
		\end{tikzpicture}
		\begin{tikzpicture} 
		\begin{axis}[ 
		cycle list name=defaultcolorlist,
		legend style={at={(1.12,1.95)}, anchor={north east}},
		xlabel=$\log_2(n)$,
		grid,
		width=5cm,
		height=4cm,
		] 
		\addplot table[x=log(n),y=log(MISE)] { 
		log(n)  log(MISE) 
		16.0  -10.386671027956352 
		17.0  -10.906323234868157 
		18.0  -12.02420612295912 
		19.0  -13.030561924022795 
		20.0  -14.133836473244045 
	}; 
	\addlegendentry{MC}

	\addplot table[x=log(n),y=log(MISE)] { 
		log(n)  log(MISE) 
		15.965784284662087  -14.624122625895069 
		17.017276025914487  -15.926972640736638 
		18.0  -16.85182456550578 
		19.019550008653873  -17.704148569387396 
		20.017276025914487  -19.060495580354043 
	}; 
	\addlegendentry{Stratif}
	\addplot table[x=log(n),y=log(MISE)] { 
		log(n)  log(MISE) 
		16.0  -16.141755010858162 
		17.0  -16.589960669586148 
		18.0  -17.93502264226188 
		19.0  -18.621379892942045 
		20.0  -19.459860334907876 
	}; 
	\addlegendentry{Sobol' +LMS}
	%
	\addplot table[x=log(n),y=log(MISE)] { 
		log(n)  log(MISE) 
		16.0  -15.615982736205153 
		17.0  -16.31525007722051 
		18.0  -17.572558128094023 
		19.0  -19.077760800136044 
		20.0  -19.794678952215136 
	}; 
	\addlegendentry{Sobol'+NUS}
	%
	%
	\addplot table[x=log(n),y=log(MISE)] { 
		log(n)  log(MISE) 
		
		16.0  -16.068445597149932 
		17.0  -16.70833866190716 
		18.0  -17.67312220237025 
		19.0  -18.652494876004223 
		20.0  -19.574632883313445 
	}; 
	\addlegendentry{Lattice+baker}
	\end{axis}
		\end{tikzpicture}
		
		\begin{tikzpicture} 
		\begin{axis}[ 
		cycle list name=defaultcolorlist,
		legend style={at={(1.12,1.95)}, anchor={north east}},
		xlabel=$\log_2(n)$,
		ylabel= $\log_2(\Var)$,
		grid,
		width=5cm,
		height=4cm,
		] 
	\addplot table[x=log(n),y=log(MISE)] { 
		log(n)  log(MISE) 
		16.0  -8.666224412165857 
		17.0  -9.416908914661668 
		18.0  -10.364919448061185 
		19.0  -11.851388003777377 
		20.0  -12.002752101335979 
		
	}; 
	\addlegendentry{MC}
	\addplot table[x=log(n),y=log(MISE)] { 
		log(n)  log(MISE) 
		16.0  -13.936914603457174 
		16.99171005377434  -15.218176303204794 
		18.094247824228052  -16.614081351574377 
		19.019550008653873  -17.489456260180614 
		20.0  -19.09614229889535 
		
	}; 
	\addlegendentry{Stratif}
	
	\addplot table[x=log(n),y=log(MISE)] { 
		log(n)  log(MISE) 
		16.0  -21.395980476569648 
		17.0  -23.219339466764758 
		18.0  -24.72212845665191 
		19.0  -26.248759822966957 
		20.0  -27.845069680521398 
		
	}; 
	\addlegendentry{Sobol'+LMS}
	%
	%
	\addplot table[x=log(n),y=log(MISE)] { 
		log(n)  log(MISE) 
		16.0  -21.868716942913437 
		17.0  -22.95962007213291 
		18.0  -24.644489563894624 
		19.0  -26.714828624742733 
		20.0  -27.315704799495244 
		
	}; 
	\addlegendentry{Sobol'+NUS}
	%
	%
	\addplot table[x=log(n),y=log(MISE)] { 
		log(n)  log(MISE) 
		
		16.0  -20.798551187185783 
		17.0  -22.763388846461144 
		18.0  -24.1031361652644 
		19.0  -25.131520772384036 
		20.0  -27.122142807184268 
		
	}; 
	\addlegendentry{Lattice+baker}
	\end{axis}
		\end{tikzpicture}
		\begin{tikzpicture} 
		\begin{axis}[ 
		cycle list name=defaultcolorlist,
		legend style={at={(1.12,1.95)}, anchor={north east}},
		xlabel=$\log_2(n)$,
		grid,
		width=5cm,
		height=4cm,
		] 
			\addplot table[x=log(n),y=log(MISE)] { 
			log(n)  log(MISE) 
			
			16.0  -8.652375791345808 
			17.0  -9.400324242995348 
			18.0  -10.442470845686135 
			19.0  -11.893633218393212 
			20.0  -11.999624023246122 
			
		}; 
		\addlegendentry{MC}
		\addplot table[x=log(n),y=log(MISE)] { 
			log(n)  log(MISE) 
			
			16.0  -13.799037473076872 
			16.99171005377434  -15.062171189680967 
			18.094247824228052  -16.42932213727655 
			19.019550008653873  -17.271389965000512 
			20.0  -18.922303261782215 
			
		}; 
		\addlegendentry{Stratif}
		
		\addplot table[x=log(n),y=log(MISE)] { 
			log(n)  log(MISE) 
			
			16.0  -20.65760855097481 
			17.0  -21.06545770985458 
			18.0  -24.03400754947537 
			19.0  -25.17007074153765 
			20.0  -26.895021907415103 
			
		}; 
		\addlegendentry{Sobol'+LMS}
		%
		%
		%
		\addplot table[x=log(n),y=log(MISE)] { 
			log(n)  log(MISE) 
			16.0  -20.449072414212182 
			17.0  -21.304884473263733 
			18.0  -23.945936869667676 
			19.0  -24.708858990566913 
			20.0  -27.32581432571954

		}; 
		\addlegendentry{Sobol'+NUS}
		%
		%
		%
		\addplot table[x=log(n),y=log(MISE)] { 
			log(n)  log(MISE) 
			
			16.0  -19.035899395738397 
			17.0  -20.943848224133713 
			18.0  -22.84080318461193 
			19.0  -24.44129308353856 
			20.0  -26.030481966065782 
			
		}; 
		\addlegendentry{Lattice+baker}
		\end{axis}
		\end{tikzpicture}
		\begin{tikzpicture} 
		\begin{axis}[ 
		cycle list name=defaultcolorlist,
		legend style={at={(1.12,1.95)}, anchor={north east}},
		xlabel=$\log_2(n)$,
		grid,
		width=5cm,
		height=4cm,
		] 
		\addplot table[x=log(n),y=log(MISE)] { 
		log(n)  log(MISE) 
		16.0  -8.53934570240081 
		17.0  -9.476899800633442 
		18.0  -10.40309577892365 
		19.0  -11.505510626641803 
		20.0  -12.178521386358495 
	}; 
	\addlegendentry{MC}

	\addplot table[x=log(n),y=log(MISE)] { 
		log(n)  log(MISE) 
		15.965784284662087  -12.80656264374196 
		17.017276025914487  -14.197559620420003 
		18.0  -15.705155699273572 
		19.019550008653873  -16.782354981277724 
		20.017276025914487  -17.855870032217776 
	}; 
	\addlegendentry{Stratif}
	\addplot table[x=log(n),y=log(MISE)] { 
		log(n)  log(MISE) 
		16.0  -16.28738135858736 
		17.0  -17.319042011109747 
		18.0  -18.883818810162555 
		19.0  -19.44211050869576 
		20.0  -20.93951862124304 
	}; 
	\addlegendentry{Sobol' +LMS}
	%
	\addplot table[x=log(n),y=log(MISE)] { 
		log(n)  log(MISE) 
		16.0  -16.117677581400883 
		17.0  -17.776532878357298 
		18.0  -18.24325663944454 
		19.0  -19.5599986210728 
		20.0  -20.53025947107574 
	}; 
	\addlegendentry{Sobol'+NUS}
	%
	%
	\addplot table[x=log(n),y=log(MISE)] { 
		log(n)  log(MISE) 
		16.0  -16.172081320286516 
		17.0  -17.111259607558328 
		18.0  -18.81875995025613 
		19.0  -19.908594307500554 
		20.0  -20.378405160721815 
		
	}; 
	\addlegendentry{Lattice+baker}
	\end{axis}
		\end{tikzpicture}
		
\caption{Plots of empirical $\log_{2}\Var[\mua]$ vs $\log_2(n)$ for various sorts and point sets, 
  based on $m=100$ independent replications, for the Heston model.
	Asian option (above) and European option (below),
	with split sort  (left),  batch sort  (middle), and Hilbert sort (right).}
\label{fig:heston-array}
\end{figure}

Following \citeN{vGIL08a}, to reduce the bias due to the discretization, we make the change of variable 
$W(t)=e^{\lambda t}(V(t)-\sigma^{2})$, with $\d W(t) = e^{\lambda t} \xi V(t)^{1/2} \d B_{2}(t)$,
and apply the Euler method to $(S,W)$ instead of $(S,V)$.
The Euler approximation scheme with step size $\delta$ applied to $W$ gives
\[
	\widetilde{W}(j\delta) = \widetilde{W}((j-1)\delta) 
	    + e^{\lambda (j-1)\delta} \xi (\widetilde{V}((j-1)\delta) \delta)^{1/2} Z_{j,2}.
\]
Rewriting it in terms of $V$ by using the reverse identity $V(t)=\sigma^{2}+e^{-\lambda t} W(t)$, 
and after some manipulations, we obtain the following discrete-time stochastic recurrence, 
which we will simulate by Array-RQMC:
\begin{eqnarray*}
	\widetilde{V}(j\delta)  &=& \max\left[0,\, \sigma^{2} + e^{-\lambda\delta} 
	 \left(\widetilde{V}((j-1)\delta)-\sigma^{2} + \xi (\tilde{V}((j-1)\delta) \delta)^{1/2} Z_{j,2}\right)\right],\\
	\widetilde{S}(j\delta)  &=& (1 + r\delta) \widetilde{S}((j-1)\delta) 
	      + (\widetilde{V}((j-1)\delta) \delta)^{1/2} \widetilde{S}((j-1)\delta) Z_{j,1},
\end{eqnarray*} 
where $(Z_{j,1},Z_{j,2})$ is a pair of standard normals with correlation $\rho$.
We generate this pair from a pair $(U_{j,1}, U_{j,2})$ of independent $\dUnif(0,1)$ variables via
$Z_{j,1} = \Phi^{-1}(U_{j,1})$ and $Z_{j,2} = \rho Z_{j,1} + \sqrt {1-\rho^{2}} \,\Phi^{-1}(U_{j,2})$. 
We then approximate each $S(j\delta)$ by $\widetilde S(j\delta)$.
The running average $\bar S_j$ at step $j$ must be the average of the $S(t_k)$ at the observation times 
$t_k \le w_j = j\delta$.  If we denote $N_j = \sum_{k=1}^c \II[t_k \le j\delta]$, we have
$\bar S_j = (1/N_j) \sum_{k=1}^{N_j} S(t_k)$, which we approximate by 
$\overline{\bar S}_j = (1/N_j) \sum_{k=1}^{N_j} \widetilde S(t_k)$.
%
Here, the state of the chain is $X_{j} = (\widetilde S(j\delta), \widetilde V(j\delta))$  
when pricing the European option and 
$X_{j} = (\widetilde S(j\delta), \widetilde V(j\delta), \overline{\bar{S}}_{j})$
when pricing the Asian option.
And two uniform random numbers, $(U_{j,1}, U_{j,2})$, are required at each step of the chain.
We thus need four-dimensional RQMC point sets for the European option and 
five-dimensional RQMC point sets for the Asian option, if we do not map the state to a 
lower-dimensional representation.
If we map the state to one dimension, as in the Hilbert curve sort, then we only need three-dimensional
RQMC points for both option types.

We ran experiments with $T=1$ (one year), $ K=100 $, $ S(0)=100 $, $ V(0)=0.04 $, $ r=0.05 $, $ \sigma=0.2 $, 
$\lambda=5 $, $ \xi=0.25 $, $ \rho=-0.5 $, and $c = \tau=16$.  
\hpierre{Maybe these numbers will change?  Perhaps $\tau=256$ and $c=12$ or 24?}%
This gives $\delta = 1/16$, 
so the time discretization  for Euler is very coarse, but a smaller $\delta$ gives similar results
in terms of variance reduction by Array-RQMC.   
\hpierre{We should report a few experiments with much smaller $h$, for example $h = 1/256$ and $h = 1/1024$,
  and compare the estimates (for expected payoff) and variance reductions.}%
\hamal{it's takes a lot of times may be I will have the results for tomorrow for  $h = 1/256$ }%
Table~\ref{tab:heston-array} reports the estimated slopes $\hat\beta$ and VRF20, 
as in Table~\ref{tab:vg-array}.
Again, we observe large variance reductions and improved convergence rates from Array-RQMC.
The best results are obtained with the split sort.
Figure~\ref{fig:heston-array} shows plots of $\log_{2}\Var[\mua]$ vs $\log_2(n)$ for selected sorts.

We tried an alternative Markov chain definition in which the chain advances by one step 
each time a uniform random number is used, as in the VG example, to reduce the dimension of the RQMC points,
but this gave no improvement.
\hpierre{Well, here I am surprised that it does not help, because the dimension of the RQMC points 
  is really reduced by 1 at all steps, if done properly!}
\hpierre{We also tried a map that combines the two state variables 
$\widetilde S(j\delta)$ and $\overline{\bar{S}}_j$ to a single variable 
$Z_j = h_j(\widetilde S(j\delta),\overline{\bar{S}}_j)$ before performing a split sort based on 
$(Z_j, \widetilde V(j\delta)))$ at step $j$,
to reduce the dimension of the sort by 1, as in the VG example, but it did not help.}
\hpierre{The VRF of MC with respect to MC should be 1. Why are we getting 1.27 and 1.28?   Very suspicious.}

\section{OPTION PRICING UNDER THE ORNSTEIN-UHLENBECK VOLATILITY MODEL}
\label{sec:ou}


	\begin{table}[hbt]
		\centering
		\small
		\caption{
		Regression slopes $\hat\beta$ for $\log_{2} \Var[\mua]$ vs $\log_2(n)$, 
		and VRF compared with MC for $n=2^{20}$, denoted VRF20, for the European and Asian  
		options under the Ornstein-Uhlenbeck model.}
		\label{tab:ou-array} 
		\begin{tabular}{l|c|c|r|c|r}
			\hline
			\multicolumn{2}{c|}{}  & \multicolumn{2}{c|}{European}& \multicolumn{2}{c}{Asian}   \\
			\hline
			Sort &  Point sets & $ \hat{\beta} $ & VRF20 & $ \hat{\beta} $ & VRF20  \\
			\hline
			\multirow{5}{2.5cm}{Batch sort 
			}	 	& MC & -1 &  1 & -1 & 1  \\
			& Stratif & -1.28 & 111 &   -1.23 &  29. \\
			& Sobol'+LMS& -1.35  &    61516   & -1.22 & 4558  \\
			& Sobol'+NUS & -1.31  &   56235   & -1.22& 5789 \\
			& Lattice+baker & -1.37  & 61318        & -1.20& 5511 \\
			\hline
			\multirow{5}{2.5cm}{Hilbert  sort  (with logistic map)
			}	 	& MC & -1 & 1  & -1 &  1  \\
			& Stratif & -1.40 &  440 &  -1.37  &  250 \\
			& Sobol'+LMS & -1.52  & 194895  & -1.40  &  41100 \\
			& Sobol'+NUS & -1.68 & 191516   & -1.37 &  39861  \\
			& Lattice+baker & -1.59  & 165351   & -1.47  &  37185  \\
			\hline
		\end{tabular} 
	\end{table}

The Ornstein-Uhlenbeck volatility model is defined by the following 
stochastic differential equations:
\begin{eqnarray*}
	\d S(t) &=& rS(t)\d t + e^{V(t)} S(t) \d B_{1}(t),\\
	\d V(t) &=& \alpha (b-V(t))\d t + \sigma \d B_{2}(t),
\end{eqnarray*}
for $ t \geq 0 $, where $ (B_{1},B_{2}) $ is a pair of standard Brownian motions with correlation $ \rho $ between  them, $ r $ is the risk-free rate, $ b $ is the long-term average volatility, 
$ \alpha $ is the rate of return to the average volatility, 
and is $ \sigma $ a variance parameter for the volatility process.
The processes $S = \{S(t),\, t\ge 0\}$ and $V = \{V(t),\, t\ge 0\}$ represent the asset price and 
the volatility process.
We simulate these processes using Euler's method with $\tau$ time steps of length $\delta$,
as we did for the Heston model, but without a change of variable.   
The discrete-time approximation of the stochastic recurrence is
\begin{eqnarray*}
	\widetilde{S}(j\delta) &=& \widetilde{S}((j-1)\delta) + r\delta \widetilde{S}((j-1)\delta)
	     + \exp\left[{\widetilde{V}((j-1)\delta)}\right] \sqrt {\delta} Z_{j,1},\\
	\widetilde{V}(j\delta) &=& \alpha \delta  b + (1-\alpha \delta  ) \widetilde{V}((j-1)\delta)
	     + \sigma  \sqrt {\delta} Z_{j,2},
\end{eqnarray*} 
where $(Z_{j,1},Z_{j,2})$ is a pair of standard normals with correlation $\rho$.
To generate this pair, we generate independent $\dUnif(0,1)$ variables $(U_{j,1}, U_{j,2})$, and put
$Z_{j,1} = \Phi^{-1}(U_{j,1})$ and $Z_{j,2} = \rho Z_{j,1} + \sqrt {1-\rho^{2}} \,\Phi^{-1}(U_{j,2})$. 
For either the European or Asian option, the state of the Markov chain and the dimension of the RQMC points
are the same as for the Heston model.
	
We ran a numerical experiment with $ T=1 $, $ K=100 $, $ S(0)=100 $, $ V(0)=0.04 $, $ r=0.05 $, 
$ b=0.4 $, $ \alpha=5 $, $ \sigma=0.2 $, $ \rho=-0.5$, and $c = \tau =16 $ (so $\delta = 1/16$).
Table~\ref{tab:ou-array} reports the estimated regression slopes $\hat\beta$ and VRF2.
\section*{CONCLUSION}
\label{sec:conclusion}

We have shown how Array-RQMC can be applied for pricing options under stochastic volatility 
models, and gave detailed examples with the VG, Heston, and Ornstein-Uhlenbeck models.
With the models, the method requires higher-dimensional RQMC points than with the simpler GBM model
studied previously, and when time has to be discretized to apply Euler's method, the number of steps
of the Markov chain is much larger.
For these reasons, it was not clear a priori if Array-RQMC would be effective.
Our empirical results show that it brings very significant variance reductions compared with 
crude Monte Carlo.

\section*{ACKNOWLEDGMENTS}

This work has been supported by a discovery grant from NSERC-Canada, a Canada Research Chair,
and a Grant from the IVADO Fundamental Research Program, to P. L'Ecuyer.

\bibliographystyle{wsc}
\bibliography{vrt,ift,fin,simul,stat,random}

\section*{AUTHOR BIOGRAPHIES}
	
\noindent {\bf AMAL BEN ABDELLAH} is a PhD student  in computer science
	at the Universit\'e de Montr\'eal, Canada. Her main research interests are randomized quasi-Monte Carlo methods,
	Array-RQMC, density estimation   and stochastic simulation  in general. She is currently
	working on quasi-Monte Carlo methods for the simulation of Markov chains and density estimation. 
	Her email address is 	\email{amal.ben.abdellah@umontreal.ca}.
\medskip

\noindent {\bf PIERRE L'ECUYER} is a Professor in the 
	Departement d'Informatique et de Recherche Op\'erationnelle, 
	at the Universit\'e de Montr\'eal, Canada. 
	He is a member of the CIRRELT and GERAD research centers. His main research interests are random
	number generation, quasi-Monte Carlo methods, efficiency improvement via variance reduction, sensitivity
	analysis and optimization of discrete-event stochastic systems, and discrete-event simulation in general. He
	has published over 270 scientific articles, 
  and has developed software libraries
	and systems for random number generation and stochastic simulation (SSJ, TestU01, RngStreams, etc.). 
	He has been a referee for 150 different scientific journals. 
More information can be found on his web page:
\href{http://www.iro.umontreal.ca/~lecuyer}%
{\url{http://www.iro.umontreal.ca/~lecuyer}}.
Email: \url{lecuyer@iro.umontreal.ca}.
\medskip

\noindent {\bf FLORIAN PUCHHAMMER} is a postdoctoral fellow  at the Universit\'e de Montr\'eal, Canada. 
	His main research interests are  discrepancy theory, quasi-Monte Carlo methods, 
	randomized quasi-Monte Carlo methods, uniform distribution of sequences, information based complexity. 
	He is currently
	working on the quasi-Monte Carlo Methods for the simulation of Markov Chains and for density estimation. 
	His email address is 	\email{florian.puchhammer@umontreal.ca}.
	
\end{document}